\documentclass[12pt]{amsart}
\usepackage{amsfonts,amsmath,amssymb}
\usepackage{graphicx}
\usepackage{tikz}
\usepackage{authblk}

\newcommand{\remove}[1]{}

\newtheorem{satz}{Theorem}[section]

\newtheorem{theorem}[satz]{Theorem}
\newtheorem{lemma}[satz]{Lemma}

\newtheorem{corollary}[satz]{Corollary}

\newtheorem{claim}[satz]{Claim}
\numberwithin{equation}{section}

\def\la{\lambda}

\def\C{\mathbb{C}}
\def\R{\mathbb{R}}

\def\({\big (}
\def\){\big )}
\def\rank{{\rm rank}}

\def\le{\leqslant}
\def\ge{\geqslant}
\def\_phi{\varphi}
\def\eps{\varepsilon}

\def\la{\lambda}

\usepackage{lipsum}
\makeatletter
\g@addto@macro{\endabstract}{\@setabstract}
\newcommand{\authorfootnotes}{\renewcommand\thefootnote{\@fnsymbol\c@footnote}}%
\makeatother

\begin{document}

\begin{center}
  \LARGE 
  Smaller Gershgorin disks for multiple eigenvalues for complex matrices\par \bigskip

  \normalsize
  \authorfootnotes
  Imre B\'ar\'any\textsuperscript{1} and J\'ozsef Solymosi\textsuperscript{2}
  \par \bigskip

  \textsuperscript{1}Alfr\'ed R\'enyi Institute of Mathematics, Budapest \par
  \textsuperscript{2}University of British Columbia, Vancouver and Obuda University, Budapest\par \bigskip

  \today
\end{center}



\begin{abstract}
Extending an earlier result for real matrices we show that multiple eigenvalues of a complex matrix lie in a reduced Gershgorin disk.  One consequence is a slightly better estimate in the real case.  Another one is a geometric application. Further results of a similar type are given for normal and almost symmetric matrices.
\end{abstract}
\section{Introduction}

Gershgorin's classic result \cite{GE} has been a major tool to estimate the eigenvalues of an $n \times n$ (complex) matrix $A=(a_{i,j})_{i,j=1}^n$ for the last 90 years. It says that if $\lambda$ is an eigenvalue of $A$, then there is $i \in [n]:=\{1,\ldots,n\}$ such that
\[
|a_{i,i}-\lambda| \le \sum_{j\ne i} |a_{i,j}|.
\]
In other words $\lambda$ lies in the {\sl Gershgorin disk} $D(a_{i,i},R_i)$ in the complex plane for some $i \in [n],$ where $a_{i,i}$ is the centre, and $R_i= \sum_{j\ne i}|a_{i,j}|$ is the radius of the disk. Equivalently, every eigenvalue $\la$ of $A$ satisfies 
\begin{equation}\label{eq:gersh}
\la \in \bigcup_{i=1}^nD(a_{i,i},R_i).
\end{equation}
In other words $\la \notin \bigcup_{i=1}^nD(a_{i,i},R_i)$ implies that $\rank(A-\la I)=n,$
where $I$ stands for the $n \times n$ unit matrix.

\medskip
In particular, if
\begin{equation}\label{eq:gersh0}
|a_{i,i}| > \sum_{j\ne i} |a_{i,j}| \mbox{ for every }i \in [n],
\end{equation}
then $0 \notin \bigcup_{i=1}^nD(a_{i,i},R_i)$ and so the matrix is not singular. Matrices satisfying (\ref{eq:gersh0}) are called 
{\em diagonally dominant} matrices. This is a classical application of Gershgorin's theorem: diagonally dominant matrices have full rank, that is, zero is not an eigenvalue because no Gershgorin disk contains the origin.  There are various generalizations of Gershgorin's theorem see for instance Alon~\cite{Alon1} and Alon and Solymosi~\cite{AS} that imply this case of diagonally dominant matrices.

\medskip
Under some special conditions, the Gershgorin bound has recently been improved by B\'ar\'any and Solymosi \cite{BS} followed  by Hall and Marsli \cite{HM}. Assume $A$ is an $n\times n$ real matrix and let $t_i$ denote the median of the numbers $a_{i,1}, \ldots, a_{i,i-1}, 0, a_{i,i+1}, \ldots, a_{i,n}$ (here $a_{i,i}$ is replaced by $0$) and set $r_i=|t_i|+\sum_{j\ne i}|a_{i,j}-t_i|.$ It is easy to see (check Section~\ref{sec:real} for an argument) that both $t_i$ and $r_i$ can be determined explicitly.

\medskip
The following result is implicit in~\cite{BS}; see the inequality in the last lines of page 3 there. A detailed proof appears in Hall and Marsli \cite{HM}. 

\begin{theorem}\label{th:BS}
Assume $\lambda$ is an eigenvalue of the real matrix $A$ whose geometric multiplicity is at least 2. Then $\la \in \bigcup_1^n D(a_{i,i},r_i).$
\end{theorem}

A similar result was proved much earlier,  in 1954,  by Ky Fan and Hoffman~\cite{FH}.  They again require geometric multiplicity,  but allow (geometric) multiplicity higher than $2$ and the corresponding radius is given in different terms. They say that it was Hadamard who discovered the basic case of Gershgorin's theorem, namely that condition (\ref{eq:gersh}) implies $\rank\, A =n$, .

\medskip
The main aim of this paper is to extend Theorem~\ref{th:BS} to complex matrices and, simultaneously, to get rid of the condition on ``geometric multiplicity''.  The proof uses Gershgorin's original theorem and is fairly simple.  It also gives a stronger result (Corollary~\ref{cor:BS} below) in the real case.  Another target is the extension of Theorem~\ref{th:BS} in the real case by replacing the smaller Gershgorin disks with slightly larger disks when the underlying matrix is normal,  and when it is almost symmetric. This is carried out in Sections~\ref{sec:normal} and \ref{sec:symm}. 

\medskip
Here comes the statement in the complex case.

\begin{theorem}\label{th:complex} Let $A$ be an $n\times n$ complex matrix and $\la$ be a complex number. If there are $c_i\in \C$ satisfying the inequality
\begin{equation}\label{eq:rank2}
|a_{i,i}-c_i-\la| >\sum_{i\neq j}|a_{i,j}-c_i| \mbox{ for all } i\in [n], 
\end{equation}
then $\rank\, (A-\la I) \ge n-1.$ 
\end{theorem}

The relation to Gershgorin disks can be seen as follows.  Fix the numbers $c_i$ satisfying (\ref{eq:rank2}) (if they exist) and define $\rho_i=\sum_{i\neq j}|a_{i,j}-c_i|$ for $i\in [n].$ Theorem~\ref{th:complex} states then, in a form analogous to (\ref{eq:gersh}), that if $\la$ is an eigenvalue of $A$ with multiplicity at least two, then $\la \in \bigcup_1^n D(a_{i,i}-c_i,\rho_i)$. Note that in this formulation there is no need for any condition on geometric multiplicity.

\smallskip
Checking whether condition (\ref{eq:rank2}) holds is fairly easy in the real case as we shall see in Section~\ref{sec:real}. The complex case is more complex and we give some comments on it at the end of Section~\ref{sec:real}. 

\smallskip
Choosing the optimal $c_i$ often reduces the radius of the Gershgorin disk for eigenvalues with a multiplicity of at least two. The $n\times n$ all one matrix, $J_n,$ shows that higher multiplicity might not guarantee a smaller radius  in Theorem \ref{th:BS} because the zero eigenvalue of $J_n$ has multiplicity $n-1$ in $J_n$ and still, it lies on the boundary of the disk with radius $r_i=1$ with center $a_{i,i}=1.$

\smallskip
Assume now that $A$ is a real $n \times n$ matrix. For all $i \in [n]$ let $t_i^*$ denote the median of the numbers $a_{i,1},\ldots,a_{i,i-1},a_{i,i+1},\ldots,a_{i,n}$ and set $r_i^*=\sum_{j\ne i}|a_{i,j}-t_i^*|.$ Theorem~\ref{th:complex} when applied to the real matrix $A$ gives the following result.

\begin{corollary}\label{cor:BS}
Assume $A$ is a real $n \times n$ matrix. If $\la \notin \bigcup_1^n D(a_{i,i},r_i^*)$, then $\rank\,(A-\la I)\ge n-1$. In other words,  if $\la$ is an eigenvalue of $A$ with multiplicity two or more, then $\la \in\bigcup_1^nD(a_{i,i},r_i^*)$.
\end{corollary}

This is a  stronger version of Theorem~\ref{th:BS} in the sense that $\la$ is a usual (algebraic) eigenvalue of $A.$ At the same time the $t_i^*$ (and then $r_i^*$) given here are slightly different from the $t_i$ in Theorem~\ref{th:BS} because  $t_i$ is the median of $n$ numbers $0,a_{i,1},\ldots,a_{i,n}$,   where $a_{i,i}$ is missing while $t_i^*$ is the median of $n-1$ numbers $a_{i,1},\ldots,a_{i,n}$ ($a_{i,i}$ is missing again).  However, Corollary~\ref{cor:BS} is stronger than Theorem~\ref{th:BS} in this sense as well because $r_i^*\le r_i$.  

\medskip
In general, bounding the rank is important in linear algebra, such bounds have various applications in different parts of mathematics. There are nice examples of applications to combinatorics in Alon's classical papers \cite{Alon2,Alon1} that provide rank bounds for real matrices where the diagonal elements are larger than other entries in their row, but not large enough for a direct application of Gershgorin's theorem.

\section{Proof of Theorem~\ref{th:complex}}\label{sec:complex}

We begin with a simple reduction. When $\la=0$ Theorem~\ref{th:complex} takes the following form.

\begin{lemma}\label{rank_bound_1}
Let  $A$ be an $n\times n$ complex matrix. Assume there are numbers $c_i\in \C$ satisfying the inequality
\begin{equation}\label{eq:rank}
|a_{i,i}-c_i| >\sum_{i\neq j}|a_{i,j}-c_i|,
\end{equation}
for all $i\in [n]$. Then $\rank\,A \ge n-1.$ \qed
\end{lemma}

The lemma implies the general case by applying it to the matrix $A-\la I.$ Indeed, the diagonal entries of this matrix are $a_{i,i}-\la$ and the off-diagonal ones are simply $a_{i,j}$ so condition (\ref{eq:rank2}) for $A$ is the same as condition (\ref{eq:rank}) for $A-\la I.$

\medskip
The {\bf proof} of the lemma is simple and uses diagonally dominant matrices.  Let $C$ be the matrix whose every entry in row $i$ is $c_i$ for $i \in [n]$ and set $D=A-C.$ Condition (\ref{eq:rank}) says exactly that $D$ is a diagonally dominant matrix. The original Gershgorin theorem shows then that $\rank\,D =n.$ As $\rank\, A\ge \rank\, D -1=n-1,$ the result follows. \qed

\medskip
We note that when condition (\ref{eq:rank2}) holds for $c_i=0$ we get back Gershgorin's circle theorem because then $C$ is the zero matrix that has rank zero.  

\medskip
{\bf Remark.} One can go one step further. Define $C$ as above (with parameters $c_i$) and let $E$ be the $n\times n$ complex matrix with every entry in its  $j$th column equal to $e_j$. We define $D=A-C-E$, so the entry $i,j$ of $D$ is $a_{i,j}-c_i-e_j$.  Again, if the parameters $c_i,e_j$ can be chosen so that $D$ is diagonally dominated, then $\rank\, A \ge d-2.$

\medskip
\section{The real case}\label{sec:real}

For the application of Theorem~\ref{th:complex} one has to decide if suitable $c_i$ values exist or not.  This is fairly easy for real matrices as explained below.

\smallskip
Assume $A$ is a real $n \times n$ matrix and let $b_{i,1}\le \ldots \le b_{i,n-1}$ be the increasing rearrangement of the $n-1$ numbers $a_{i,1},\ldots,a_{i,n}$ where $a_{i,i}$ is missing.  The function $f(t)= \sum_{i\neq j}|a_{i,j}-t|$ is piecewise linear and convex. It attains its minimum at the median of the numbers $b_{i,j}$.  This median is the single point $t_i^*:=b_{i,n/2}$ when $n$ is even and can be taken for any point of the interval $J:=[b_{i,(n-1)/2},b_{i,(n+1)/2}]$ when $n$ is odd.  We remark that in this case, $r_i^*$ is the same no matter what $t\in J$ is chosen for the median.

\begin{claim}\label{cl:real} In the real case condition (\ref{eq:rank2}) is satisfied by some $c_i$ if and only if it is satisfied by $c_i=t_i^*$ when $n$ is even and by one of the endpoints of $J$ when $n$ is odd. 
\end{claim}

{\bf Proof.} One direction is easy: there is nothing to check when (\ref{eq:rank2}) is satisfied at $c_i=t_i^*$ or when it is satisfied by one of the endpoints of $J$.

\begin{figure}[h!]
\centering
\includegraphics[scale=0.85]{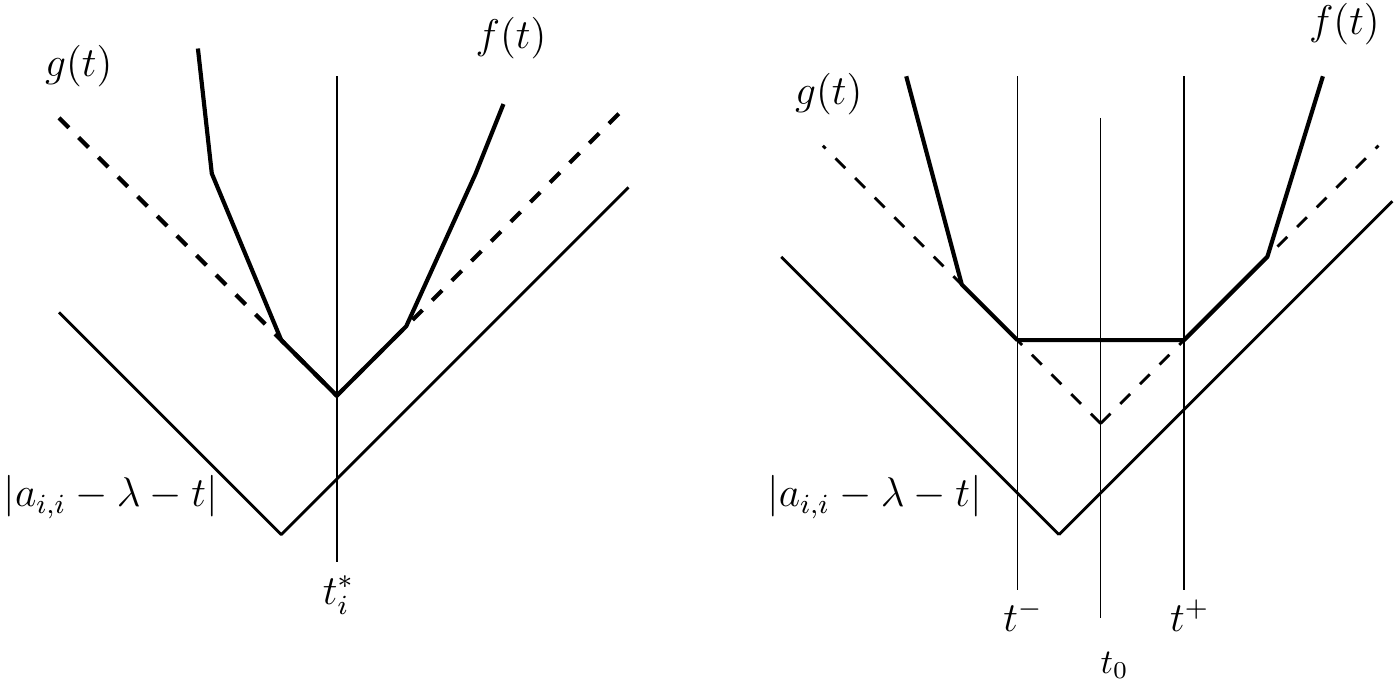}
\caption{The functions $f(t),g(t)$ and $|a_{i,i}-\la-t|$.}
\label{fig:slope}
\end{figure}

For the other direction assume first that $n$ is even and that condition (\ref{eq:rank2}) is not satisfied at $t_i^*$, that is $|a_{i,i}-t_i^*-\la| \le f(t_i^*)$.  The function $g(t)=|t-t_i^*|+f(t_i^*)$ is a translated copy of the function $t \to |t|$ and $|a_{i,i}-t-\la| \le g(t)$ for all $t$ follows from $|a_{i,i}-t_i^*-\la|\le f(t_i^*)$. Moreover $g(t)$ has slope $+1$ for $t>t_i^*$ and slope $-1$ for $t<t_i^*$, while the slope of $f(t)$ is at least $1$ for $t>t_i^*$ and and at most $-1$ for $t<t_i^*$; see Figure~\ref{fig:slope} left. Then $g(t) \le f(t)$ for all $t \in \R$ because $g(t_i^*)=f(t_i^*)$ and $|a_{i,i}-\la-t|\le f(t)$ for all $t$ follows. 

\smallskip
The case of odd $n$ is similar.  For simpler notation set $J=[t^-,t^+]$ and $t_0=\frac 12(t^-+t^+)$. Suppose that (\ref{eq:rank2}) is not satisfied at $t=t^-$ and at $t=t^+$.  Observe that $f(t^-)=f(t^+).$ Setting $g(t)=|t-t_0|+f(t^+)-\frac12 (t^+-t^-)$ we see that $g(t^-)=f(t^-)$ and $g(t^+)=f(t^+)$,  see Figure~\ref{fig:slope} right.  Then $|a_{i,i}-\la-t|\le g(t)$ for all $t$ because this inequality holds for $t=t^+$ and $t=t^-$.  Next one shows that $g(t)\le f(t)$ for all $t\in \R$. This holds because for $t>t^+$ the slope of $g(t)$ is $1$ and the slope of $f(t)$ is at least one  and for $t<t^-$ the slope of $g(t)$ is $-1$ and that of $f(t)$ is at most $-1$. Moreover $f(t)$ is the constant $f(t^+)=f(t^-)$ for $t\in J$, and there $g(t)$ is below this constant.\qed

\medskip
So in the real case in order to check whether condition (\ref{eq:rank2}) holds and for finding the suitable $c_i$s one has to compare the minimum of $f(t)$ with the value of $|a_{i,i}-\la-t|$ at one or two well-defined points.  

\medskip
Finding $c_i$ in the complex case is a different task. We have to decide if the function $t \to  |a_{i,i}-t|$ is below the function $t \to \sum_{i\neq j}|a_{i,j}-t|$ for all complex numbers $t$ or not.  One can check if this holds when $t=a_{i,j}$ for all $j\in [n]$. At every other $t \in \C$ both functions are differentiable so one could, in principle, decide if the maximum of the function
$$|a_{i,i}-t| -\sum_{i\neq j}|a_{i,j}-t|$$
is positive for any fixed $i \in [n]$ or not. As the real case indicates, a good candidate for $c_i$ is where the function $f(t)=\sum_{j\ne i}|a_{i,j}-t|$ attains its minimum on $t \in \C$. This point can be determined by convex programming.  We remark without proof that $\arg \min f(t)$ is a single point unless all the points $a_{i,j}$ ($j\ne i$) are collinear.  The latter case is covered by Claim~\ref{cl:real}.

\section{An example}\label{sec:example}

We show a geometric application of Lemma \ref{rank_bound_1} which is similar to an example of Alon  \cite{Alon1} and of Bukh and Cox~\cite{BC}. The unit $\ell_1$-ball in $\mathbb{R}^d$ is the convex hull of the $2d$ points $\pm e_i,$ i.e.,  the collection of $x\in \mathbb{R}^d$ such that $\Vert x\Vert_1\leq 1.$ Here $e_i$ denotes the vector with a 1 in the $i$th coordinate and 0's elsewhere. The Euclidean distance between $e_i$ and $-e_i$ is $2$ (for all $i$), and the distance between two other vertices is exactly$\sqrt{2}.$ What happens if we relax the distance constraints a bit? What is the maximum number of point pairs, $(p_i,q_i)$ in $\mathbb{R}^d,$ to be denoted by $k,$ such that  $(p_i-q_i)^2=4$ for every $i \in [k]$ and for every pair $i,j \in [k],\; i\ne j$ we have
\[
|(p_i-q_j)^2-{2}| \leq\varepsilon,\; |(p_i-p_j)^2-{2}| \leq\varepsilon, \; |(q_i-q_j)^2-{2}| \leq\varepsilon.
\]
The answer is that there are at most two extra point pairs if $\eps$ is small enough.  More precisely 

\begin{claim}\label{Octahedron}
If $\varepsilon< 2/(3d+5)$, then $k\leq d+2.$
\end{claim}
{\bf Proof.} Assume that $k>d+2$ and consider the first $d+3$ pairs $(p_i,q_i).$ The $(d+2)\times(d+2)$ matrix $M$ is defined by the dot products $m_{i,j}=(p_i-p_1)\cdot(q_j-p_1),$ for all $2\leq i,j\leq d+3.$
The rank of $M$ is at most $d.$ The identity $(p_i-p_1)^2+(q_j-p_1)^2-(p_i-q_j)^2=2(p_i-p_1)\cdot(q_j-p_1)$ (the cosine theorem in trigonometry) and the conditions imply that
\[
(2-\eps)+(2-\eps)-(2+\eps) \le 2(p_i-p_1)\cdot(q_j-p_1) \le (2+\eps)+(2+\eps)-(2-\eps).
\]
Consequently
\[
|m_{i,j}-1|\le \frac 32 \eps \mbox{ for } i\ne j \mbox{ and }|m_{i,i}-1|\ge 1- \eps \mbox{ for } 2\le i\le d+3.
\]
We check that Lemma~\ref{rank_bound_1} applies now with $c_i=1$ for all $i=2,3,\ldots,d+3$. Indeed $|m_{i,i}-c_i|\ge 1- \eps$ and $\sum_{j\ne i}|m_{i,j}-c_i|\le (d+1)\frac 32 \eps$ when $i\ne j.$ Thus the inequality  $|m_{i,i}-c_i|>\sum_{j\ne i}|m_{i,j}-c_i|$ follows from the condition $\varepsilon< 2/(3d+5)$ via a simple computation. By Lemma~\ref{rank_bound_1}, $\rank\,M\ge d+2-1=d+1,$ contradicting the fact that $\rank\,M\le d.$\qed

\medskip
The following construction shows that the result in Claim \ref{Octahedron} is close to being sharp. Let us choose a $d$ such that there is a Hadamard matrix, $H_n$ of order $n=d+2.$ The row vectors of $H_n$ and $-H_n$ form a scaled and rotated $\ell_1$ ball in $\mathbb{R}^n.$ The pointset we are going to consider in $\mathbb{R}^d$ consists of the projection of the $2n$ row vectors $\pm\frac{1}{\sqrt{d}}H_n$ by deleting the last two coordinates of every row. The squared distance between non-antipodal points is between $2-4/d$ and $2+4/d$ while the distance between the $d+2$ antipodal points is exactly 2.

\medskip
One can relax the condition $(p_i-q_i)^2=4$ to $|(p_i-q_i)^2-4| \leq \eps$. In this case, the proof goes along the same steps as above and gives that for $\eps <2/(3d+6)$ the maximal number of such $p_i,q_i$ pairs is at most $d+2$. We omit the details.

\section{Normal Matrices}\label{sec:normal}

Here we show that if two eigenvalues of a real and normal matrix are close, then both of them lie in a smaller Gershgorin disk. A matrix is normal if its eigenvectors belonging to distinct eigenvalues are orthogonal. For instance, a symmetric and real matrix is always normal.

\begin{theorem}\label{th:norm}
Assume $\la$ and $\mu$ are distinct eigenvalues of a real and normal $n \times n$ matrix $A$, $n\ge 3$. Then both $\la$ and $\mu$ lie in the disk $D(a_{i,i},\rho_i)$  for some $i\in [n]$ where
$\rho_i= r_i+\sqrt n|\la-\mu|$.
\end{theorem}

{\bf Proof.} Let $A$ be a normal, real, $n\times n,$ matrix. Let $v=(v_1,\ldots,v_n)$ resp. $w=(w_1,\ldots,w_n)$ be the eigenvectors corresponding to $\la$ and $\mu$, the coordinates $v_i,w_j \in \C$. As $n\ge 3$ there are  $\alpha, \beta \in \mathbb{C},$ such that
\[\sum_{i=1}^n\alpha v_i+\beta w_i=0,\]
and the largest coordinate is one, i.e. $\alpha v_i+\beta w_i=1$, and $|\alpha v_j+\beta w_j|\leq 1$ for any $j \in [n].$

Let $t=\alpha v+\beta w$. Then $A t= \mu\alpha v+\lambda\beta w$ and
\[
a_{i,1}t_1+a_{i,2}t_2\ldots +a_{i,i}+\ldots +a_{i,n}t_n=\mu\alpha {v}_i+\lambda\beta{w}_i= \mu + (\lambda-\mu)\beta{w}_i
\]
because $t_i=1$. Set $S=a_{i,1}t_1+a_{i,2}t_2\ldots +a_{i,i}+\ldots +a_{i,n}t_n-a_{i,i}$. Then
\[
|\mu -a_{i,i}| \le |S|+|\la-\mu| |\beta w_i|.
\]
Bounding $S$ comes from Theorem 2 in \cite{BS} and, in a slightly more general form, from Lemma 2.4 of \cite{HM} in the form of the following fact.
\begin{lemma}\label{l:basic} $|S| \le r_i$.
\end{lemma}
The proof is the same as in \cite{BS} and we present it at the end of this section.

\medskip
We bound the error term $|\beta w_i|$ using that $A$ is normal. The two eigenvectors are orthogonal, so we have
\[
n\geq |t|^2=|\alpha v|^2+|\beta w|^2,
\]
which implies that $|\beta w_i|\leq \sqrt{n}.$\qed

\medskip
The {\bf proof} of Lemma~\ref{l:basic} is simple.
\begin{eqnarray*}
S &=& a_{i,1}t_1+\ldots +a_{i,i-1}t_{i-1}+0t_i+a_{i,i+1}t_{i+1}+ \ldots +a_{i,n}t_n\\
   &=& (a_{i,1}-x)t_1+\ldots +(a_{i,i-1}-x)t_{i-1}+\\
   &&+(0-x)t_i+(a_{i,i+1}-x)t_{i+1}+ \ldots +(a_{i,n}-x)t_n
\end{eqnarray*}
for every $x \in \R$ because $\sum t_i=0$. As $|t_i| \le 1$ for all $i \in [n]$ this implies  that
\[
|S| \le |a_{i,1}-x|+\ldots +|a_{i,i-1}-x|+|0-x|+|a_{i,i+1}-x|+ \ldots +|a_{i,n}-x|.
\]
The function $x \to |a_{i,1}-x|+ \ldots +|a_{i,n}-x|$ is the sum of convex functions and, as we have seen, its minimum is reached on the median, $t_i$,  of the numbers
$a_{i,1},\ldots,a_{i,i-1},0,a_{i,i+1},\ldots, a_{i,n}$ and equals $r_i.$\qed

\section{Almost Symmetric Matrices}\label{sec:symm}

This section extends the previous results to almost symmetric real matrices.

There are different ways to define almost symmetric matrices. They appear in various contexts in applied linear algebra. For example in \cite{BB} and \cite{ZCL} a nearly symmetric matrix is defined as a matrix in which the overwhelming majority of entries are symmetric about its diagonal. In other examples, an almost symmetric matrix $M$ is given as $M=S+E$ where $S$ is a symmetric matrix and $E$ is a small error or ``noise''. We can't repeat the arguments on Theorem~\ref{th:norm} because even small changes can make the eigenvectors far from being orthogonal, like in the example below.

\[M=\left[
\begin{array}{cc}
 1 & \eps   \\
 0 & 1+\eps  \\
\end{array}
\right]\]

We will not give a formal definition of what an almost symmetric matrix is. Instead we introduce a parameter, $\Delta(A),$ measuring the symmetry of a matrix $A$:
\[
\Delta(A)=\max_{i \in [n]}\left|\sum_{j=1}^n(a_{i,j}-a_{j,i})\right|.
\]
In this definition the pairs $a_{i,j}$ and $a_{j,i}$ can be far from each other, the parameter $\Delta(A)$ is the maximum difference between the sum of the entries in row $i$ and that of the entries in column $i$. For instance $\Delta(A)=0$ for a doubly stochastic matrix $A$.

With this new matrix parameter, we generalize Theorem~\ref{th:norm} to arbitrary real matrices. We are going to use the well-known fact that an eigenvector of a real matrix $A$ which belongs to the eigenvalue $\lambda$ is orthogonal to any eigenvector of $A^T$ which belongs to a different eigenvalue $\mu.$ 

\remove{We include the simple argument for the sake of completeness.DELETE????
\[
vA\circ w=vAw^T=(vAw^T)^T=wA^Tv^T=wA^T\circ v,
\]
so
\[
\lambda v\circ w=\mu w\circ v, \quad (\lambda-\mu)v\circ w=0.
\]
}
\medskip

\begin{theorem}\label{th:symm} Assume that $A=\{a_{i,j}\}_{i,j=1}^n$ is a real matrix and $n\ge 3$.
If $\lambda$ and $\mu$ are two distinct eigenvalues of $A,$ then there is $i \in [n]$ such that both $\la$ and $\mu$ lie in the disk $D(a_{i,i},\rho_i)$ where $\rho_i=r_i+\sqrt n(\Delta +|\la-\mu|)$ where $r_i$ is the same as before.
\end{theorem}

{\bf Proof.} Let $v$ be the eigenvector of $\lambda$ in $A$ and $w$ be the eigenvector of $\mu$ in $A^T.$ The two vectors are orthogonal, so they are linearly independent.
Then there are $\alpha,\beta \in \C$ such that the coordinates of $\alpha v+\beta w$ add up to zero:
\[
\sum_{i=1}^n \alpha u_i +\beta w_i=0,
\]
and the largest norm coordinate (indexed by $i \in [n]$) is one. We assume again that $\alpha v_i +\beta w_i=1$
and so $|\alpha v_j +\beta w_j|\le 1$ for all $j \in [n].$ With these notations we have the system of equations
\[
A\alpha u + A^T \beta w =\la \alpha v+\mu \beta w= \lambda (\alpha v + \beta w) + (\mu - \lambda)\beta w,
\]


Let us consider the $i$th row
\[
\alpha v_1a_{1,i}+\beta w_1a_{i,1}+\ldots +a_{i,i}+\ldots  +\alpha v_ja_{j,i}+\beta w_ja_{i,j}+\ldots = \lambda + (\mu - \lambda)\beta w_i.
\]
From this we have
\[
|\lambda - a_{i,i}|\  \leq\left |\sum_{j\ne i}^n(\alpha v_j+\beta w_j)a_{j,i} - \beta\sum_{j=1}^n w_j(a_{j,i}-a_{i,j}) - (\mu - \lambda)\beta w_i \right |
\]


Using the bound $|\beta w_i|\leq \sqrt{n}$ from the previous proof we have

\begin{eqnarray*}
\begin{split}
|\lambda - a_{i,i}|\ & \leq\left |\sum_{j\ne i}(\alpha u_j+\beta w_j)a_{j,i} - \beta\sum_{j=1}^n w_j(a_{j,i}-a_{i,j}) - (\mu - \lambda)\beta w_i \right |\leq\\
& \leq \left |\sum_{j\ne i}(\alpha u_j+\beta w_j)a_{j,i}\right| + \sqrt{n}\left|\sum_{j=1}^n (a_{j,i}-a_{i,j})\right| + \sqrt{n}|(\mu - \lambda) |\leq \\
& \leq r_i + \sqrt{n}\left( \Delta(A) + \left|\mu - \lambda \right |  \right).
\end{split}
\end{eqnarray*}
Here the estimate $\left|\sum_{j\ne i}(\alpha u_j+\beta w_j)a_{j,i}\right| \le r_i$ is the same as in the proof of Theorem~\ref{th:norm}.\qed


\section{Acknowledgements}
Research of both authors was supported by a Hungarian National Research Grant KKP grant no. 133819.
Research of IB was supported in part by Hungarian National Research grants no. 131529, 131696.
Research of JS was supported in part by an NSERC Discovery grant and OTKA K  grant no.119528.

\end{document}